%% file: AnticommutingPairsArxiv.tex
\documentclass[12pt]{article}
\usepackage{verbatim}
\usepackage{amssymb}
\usepackage{amsthm}
\usepackage{amsmath}
\usepackage{mathrsfs}
\usepackage[margin=1in]{geometry}
\usepackage{mathtools}

\input def_latex5

\usepackage[usenames,dvipsnames] {color}

\definecolor{dark_purple}{rgb}{0.4, 0.0, 0.4}
\definecolor{dark_green}{rgb}{0.0, 0.7, 0.0}

\numberwithin{equation}{section}

\title{Random anti-commuting Hermitian matrices}
\author{
John E. M\raise.5ex\hbox{c}Carthy
\thanks{Partially supported by National Science Foundation Grant  
DMS 2054199}\\
Washington University in St. Louis
\and
Hazel T. M\raise.5ex\hbox{c}Carthy\\
University of Illinois at Urbana-Champaign
\\
}
\date{\today}

\newcommand{\mc}{M\raise.45ex\hbox{c}Carthy}

\usepackage{mathtools}

\newcommand\Z{\mathbb Z}
\newcommand\sn{\Sigma_n}

\newcommand\mn{{\mathbb M}_n}
\newcommand\mnd{{\mathbb M}_n^d}

\renewcommand\l{\lambda}

\newcommand\tw{{\tilde{w}}}
\newcommand\un{{\mathcal U}(n)}

\newcommand\an{{{\frak A}_n}}
\newcommand\ang{{{\frak A}_{n,{\rm gen}}}}
\newcommand\cnd{{{\frak C}_n^d}}
\newcommand\bn{{\frak B}_n}
\newcommand\rtpg{ (\R^2_+)^p_{\rm gen}}
\newcommand{\bigzero}{\mbox{\normalfont\Large 0}}
\renewcommand\ii{{\iota}}

\begin{document}

\bibliographystyle{plain}

\maketitle
\begin{abstract}
We consider pairs of anti-commuting $2p$-by-$2p$ Hermitian matrices that are chosen randomly with respect to a Gaussian measure. 
Generically such a pair decomposes into the direct sum of $2$-by-$2$ blocks on which
the first matrix has eigenvalues $\pm x_j$ and the second has eigenvalues $\pm y_j$. We call $\{ (x_, y_j) \}$ the skew spectrum of the pair. We derive a formula for the probability density of the skew spectrum, and show that the elements are repelling.

\renewcommand{\thefootnote}{\fnsymbol{footnote}} 
\renewcommand{\thefootnote}{\arabic{footnote}} 
\end{abstract}
\section{Introduction}
\label{secin}

The study of random matrices goes back at least to the 1920's, but it came to prominence in physics with the work of
 of Wigner \cite{wi55,wi57,wi58} and  Dyson \cite{dy62a,dy62b, dy62c},
 who used results from random matrices to predict the eigenvalues of complicated Hamiltonians.
 See eg. \cite{meh04} for an account.
What happens if we choose multiple Hamiltonians whose interaction forces them  to satisfy certain algebraic relations?
In \cite{mcc23a}, this question was studied when the Hamiltonians commute (see Section \ref{secb}).
The purpose of this note is to study the eigenvalue distribution of random pairs of anti-commuting Hermitian matrices.

First, let us define what we mean by a random $d$-tuple of matrices satisfying given algebraic relations.
Let $\mn$ denote the algebra of $n$-by-$n$ complex matrices, and 
let $\sn$ denote the Hermitian matrices in $\mn$. 
Let ${\frak V}_n \subseteq \mnd$ be an algebraic set, by which we mean there are non-commutative polynomials
$p_1, \dots, p_N$ in the $2d$ variables  $\{x^1, (x^1)^*, \dots, x^d, (x^d)^*\}$ so that
\be
\label{eqa01}
{\frak V}_n \= \{ X \in \mnd : p_j(X) = 0 \ \forall\ 1 \leq j \leq N \} .
\ee
The set ${\frak V}_n$ can be thought of as a subset of $\C^{dn^2} = \R^{2dn^2}$, and as such there is a natural measure on it, consisting of Hausdorff measure of the real dimension of ${\frak V}_n$. We will write this measure as $dX$. To convert this infinite measure to a probability measure, we multiply by something like a Gaussian weight.

For $X = (X^1, \dots, X^d)$ in $\mnd$ define its Frobenius (or Hilbert-Schmidt) norm by 
\[
\| X \|^2_F \ := \ \sum_{r=1}^d \sum_{i,j=1}^n  |X_{ij}^r|^2 .
\]
Let $w$ be a continuous  non-negative function on $[0,\i)$, which has enough moments that
 $w(\| X \|_F) dX$ is a finite measure on ${\frak V}_n$.  We will assume $w$ is normalized
so  $w(\| X \|_F) dX$ is a probability measure. 

\begin{definition}
A random $d$-tuple in ${\frak V}_n$ is a random variable with values in ${\frak V}_n$ and with distribution $w(\| X \|_F) dX$.
\end{definition}

In particular, in this note we will study the set
\[
\an \ := \ \{ (X,Y) \in \sn^2: XY + YX = 0 \} .
\]
We shall assume that $n = 2p$ is even.
In Section \ref{secc}, we will see that 
 $\an$ has dimension $n^2 + \frac n2$, and that generically elements of $\an$ are unitarily equivalent to a pair of the form
 \be
 \label{eqa1}
 X =
 \left(\begin{array}{@{}c|c@{}|c}
  \begin{matrix}
  x_1 & 0 \\
  0 & -x_1
  \end{matrix}
  & \bigzero& \\
\hline
  \bigzero &
  \begin{matrix}
  x_2 & 0 \\
  0 & -x_2
  \end{matrix} & \\
\hline
&&\ddots
\end{array}\right), \quad
 Y =
 \left(\begin{array}{@{}c|c@{}|c}
  \begin{matrix}
  0 & y_1 \\
  y_1 & 0
  \end{matrix}
  & \bigzero& \\
\hline
  \bigzero &
  \begin{matrix}
   0 & y_2 \\
  y_2 & 0
  \end{matrix} & \\
\hline
&&\ddots
\end{array}\right) ,
 \ee
where each $x_j$ and $y_j$ is positive. We shall call the pairs $\{ (x_j, y_j) : 1 \leq j \leq p \}$ the {\em skew spectrum} of $(X,Y)$.

Here is our main result.
\bt
\label{thma1}
Let $Z =  (X,Y)$ be chosen randomly in $\an$ with distribution $w( \| Z\|_F)$. Then the probability 
distribution of the skew spectrum of $(X,Y)$ on $\R^{2p}_+$  is given by
\begin{multline}
\label{eqa2}
\rho_n (x_1, y_1, \dots, x_p, y_p) \= \\
C_n \ w(\| Z |_F^2 ) \ \prod_{1 \leq k \leq p} \,  x_k y_k\sqrt{x_k^2 + y_k^2} \prod_{1 \leq i < j \leq p}
\left[(x_i - x_j)^2 +(y_i-y_j)^2 \right]
\left[(x_i + x_j)^2 +(y_i-y_j)^2 \right] \\
\times
\left[(x_i - x_j)^2 +(y_i + y_j)^2 \right]
\left[(x_i + x_j)^2 +(y_i +y_j)^2 \right].
\end{multline}
\et
If $x_i, x_j, y_i, y_j$ are bounded and bounded away from zero,
the last  factor in \eqref{eqa2} 
 is bounded above and below by
 \be
 \label{eqa3}
 \notag
 \left[ (x_i - x_j)^2 + (y_i - y_j)^2 \right]
 \ee
 (see Lemma \ref{lemd1}).
 This quadratic vanishing is of the same order as in the Ginibre formula for commuting Hermitian pairs, showing that the repulsion between the elements of the skew spectrum is similar to, though more complicated than, the repulsion between the joint eigenvalues for a commuting Hermitian pair.

\section{Random commuting matrices}
\label{secb}
 
 In this section, we give some results about random commuting Hermitian matrices.
 We shall not use these explicitly in the following sections, but they serve as 
 a guide to what we would like to achieve in the anti-commuting case.
 When $d=1$, Ginibre \cite{gi65} proved that for the Gaussian Hermitian ensemble,
 the eigenvalues of a random Hermitian matrix in $\sn$ have the distribution
 \be
 \label{eqb1}
 \rho(\l_1, \dots, \l_n) \= C_n  e^{- \frac 12 \sum_{j=1}^n \l_j^2} \prod_{1 \leq i < j \leq n} |\l_i - \l_j|^2 
 \ee
 on $\R^n$. We use $C_n$ to denote a constant that depends on $n$, and may vary from one
 occurrence to another. An analogous formula  to \eqref{eqb1} turns out to hold not just in the Gaussian case, but if the matrices
 are chosen with respect to any weight that depends only on $\|X\|_F$---see e.g. \cite{ta12} for an account.
Wigner proved in \cite{wi58}, subject to all moments having bounds independent of $n$, that 
if $X_n$ is chosen in $\sn$ with distribution $w(X)$ that depends only on $\|X\|_F$, then the density
of eigenvalues of  $\frac{1}{\sqrt{n}} X_n$ converges almost surely to the semi-circular distribution
\[
\frac{1}{2\pi} \sqrt{4-x^2} dx 
\]
on $[-2,2]$.
 
 Now let $d >1$, and let $\cnd$  denote the set of commuting $d$-tuples in $\sn^d$.
 Let
  $w(X)$ be a weight on $\cnd$  that depends only on $\|X \|_F$ and is normalized to have
 $\int_\cnd w(X) dX$. 
 An eigenvalue of $X$ is now a $d$-tuple in $\R^d$ (since the matrices commute, they have common eigenvectors). If $(\l_1, \dots, \l_n)$ are the eigenvalues of $X$, then 
 \[
 \sum_{j=1}^n |\l_j|^2 \= \| X \|_F^2 ,
 \]
 where $|\l | = \sqrt{\sum_{r=1}^d |\l^r|^2}$ is the Euclidean norm in $\R^d$.
 Therefore there is a function $\tw: (\R^d)^n \to \R$ so that
 \[
 w(X) \= \tw(\l) .
 \]
 
 In \cite{mcc23a} it was shown that the Ginibre formula still holds.
 \bt
 \label{thmb1}
  For $X$ in $\cnd$  with distribution $w$ as above, the eigenvalues of $X$ have density
 \be
 \label{eqb2}
 \kappa_n(\l_1, \dots, \l_n) \= C_n \, \tw(\l)  \prod_{1 \leq i < j \leq n} |\l_i - \l_j |^2 .
 \ee
 \et
 Any $X$ in $\cnd$ is unitarily equivalent to a $d$-tuple of diagonal matrices. The unitary implementing this
 is generically unique up to multiplication by a diagonal unitary. Let $\un$ denote the unitary group in $\mn$,
 and let $\T^n$ be the subgroup of diagonal unitaries. Let $\nu$ be volume measure on the homogeneous space $\un / \T^n$. Then \eqref{eqb2} asserts that the measure $w(X) dX$ decomposes as
 \be
 \label{eqb3}
 \notag
 w(X) dX \= C_n \, \tw(\l) \prod_{1 \leq i < j \leq n} |\l_i - \l_j |^2\  d\l d\nu .
 \ee
 
 Let us now restrict to the Gaussian case $w(X) = C_n e^{- \gamma \| X\|_F^2}$.
  The equilibrium measure with respect to the logarithmic potential is the probability measure $\mu_d$ that minimizes the
 logarithmic energy
 \be
 \label{eqb4}
 \notag
 I(\mu) \= \int_{\R^d} \int_{\R^d} \log \frac{1}{|x-y|} d\mu (x) d\mu(y) + \int_{\R^d} \gamma |x|^2 d \mu (x) .
\ee
The equilibrium measure exists, is unique, and is compactly supported \cite[Thm. 4.4.14]{bhs19}.
The eigenvalue density, scaled by $\frac{1}{\sqrt{n}}$, converges to the equilibrium measure.
We shall let ${\mathbb E}_n$ denote expectation at the $n^{\rm th}$ level of the process.
\bt
\label{thmb2} \cite{mcc23a}
Let  $X_n$ be chosen in $\cnd$ with distribution $C_n e^{-\gamma \| X \|^2_F} dX$.
Let $\phi$ be a continuous bounded function on $\R^d$. Then
\be
\label{eqb5}
\notag
\lim_{n \to \i} {\mathbb E}_n \left[ \frac 1n {\rm tr} (\phi(\frac{1}{\sqrt{n}}X_n) )\right]
\= \int_{\R^d} \phi(x)  d\mu_d(x) .
\ee
\et

In this Gaussian case, the equilibrium measures have been calculated explicitly by Chafai, Saff and Womersley \cite{csw22,csw23}. We let $\sigma^{d-1}_{R_d}$ denote normalized surface area on the sphere of radius $R_d$ in $\R^d$.
 \bt
 \label{thmb3}
 [Chafai, Saff, Womersley]
 Let $Q(x) = \gamma |x|^2$. Then the equilibrium measure is supported on the ball of radius $R_d$, and  is given by
\begin{eqnarray}
\notag
\label{eqb6}
 \frac{2}{\pi R_1^2} \sqrt{(R_1^2 - x^2)_+}\ dx, &\quad R_1 := \sqrt{\frac 2\gamma}, & d = 1;
 \\
 \frac{1}{\pi R_2^2} 1_{|x| < R_2}\  dx^1 dx^2,
 \label{eqb7}
 \notag
 &R_2 := \frac{1}{ \sqrt{\gamma}}, & d =2;
 \\
 \frac{1}{\pi^2 R_3^2} \frac{1}{\sqrt{R_3^2 - |x|^2}} 1_{|x| < R_3}\ dr d\sigma^2_1,
 &R_3 := \sqrt{\frac{2}{3\gamma}},& d = 3;
 \label{eqb8}
 \notag
 \\ \sigma^{d-1}_{R_d},
 & R_d := \frac{1}{\sqrt{2\gamma}},&d \geq 4.
 \label{eqb9}
 \notag
 \end{eqnarray}
 \et
 
 \section{Generic elements in $\an$}
 \label{secc}

 Let $\bn$ denote $\{ (X,Y) \in \mn^2 : XY +YX = 0 \}$. 
 The set of commuting pairs in $\mn$ is an irreducible variety \cite{mt55}, but $\bn$ is not.
  In \cite{cw20} Chen and Wang showed that for each triple $(q,m,r)$ of non-negative integers that satisfy
  \[ 2q + m + r \= n \]
  there is an irreducible variety ${\frak Z}_{q,m,r}$ so that
  \[
  \bn \= \bigcup_{ 2q+m + r = n}{\frak Z}_{q,m,r} .
  \]
 These varieties can be described as follows.
 Let ${\frak U}_{q,m,r}$ be the set of pairs $(X,Y)$ in $\bn$ that are jointly similar to a block-diagonal pair of the form
 
  \begin{eqnarray}
  \notag
 \label{eqc1}
 X &=&
 \left(\begin{array}{@{}c|c|c@{}|c|c}
  \begin{matrix}
  x_1 & 0 \\
  0 & -x_1
  \end{matrix}
  & && &\\
  \hline & \ddots & &&\\
\hline
  & &
  \begin{matrix}
  x_q & 0 \\
  0 & -x_q
  \end{matrix} && \\
\hline
&&&A_m& \\
\hline
&&&&\bigzero_r
\end{array}\right), \\
  Y &=&
 \left(\begin{array}{@{}c|c|c@{}|c|c}
  \begin{matrix}
 0 & y_1 \\
  z_1 & 0
  \end{matrix}
  & && &\\
  \hline & \ddots & &&\\
\hline
  & &
  \begin{matrix}
 0 & y_q \\
  z_q & 0
  \end{matrix} && \\
\hline
&&&\bigzero_m& \\
\hline
&&&& B_r
\end{array}\right),
 \label{eqc2}
 \notag
 \end{eqnarray}
 where $A_m$ and $B_r$ are arbitrary diagonal matrices, of size $m$-by-$m$ and $r$-by-$r$ respectively,
 $X$ has rank $2q+m$, $Y$ has rank $2q+r$,  all the non-zero eigenvalues of $X$ are distinct, and all the non-zero eigenvalues of $Y$ are distinct.
 Then ${\frak Z}_{q,m,r}$ equals the Zariski closure of ${\frak U}_{q,m,r}$, and has complex dimension $n^2 + q$ \cite{cw20}.

 What does this tell us about $\an$? 
 As $\an = \bn \cap \sn^2$, we have
 \be
 \notag
 \label{eqc3}
 \an \= \bigcup_{2q+m+r =n} {\frak Z}_{q,m,r} \cap \sn^2 .
 \ee
 Similar arguments to the ones given in \cite{cw20} show that the real dimension of $ {\frak Z}_{q,m,r} \cap \sn^2$ is $n^2 +q$, so only the largest component,
 $ {\frak Z}_{p,0,0} \cap \sn^2$ will have positive measure with respect to $w(\| X \|_F) dX$. So we can restrict our
 attention to what we will call the generic elements in $\an$, namely the set of full measure consisting of pairs 
 that are jointly unitarily equivalent to some $(X,Y)$ as in \eqref{eqa1} with $\{ x_1, \dots, x_p \}$ and
 $\{ y_1, \dots , y_p \}$ both consisting of $p$ distinct positive numbers. (Since $Y$ is self-adjoint, we have
 $\bar z_j = y_j$, and we can choose $y_j$ to be positive by conjugating the $j^{\rm th}$ block by an appropriate diagonal unitary.)
 The skew spectrum of such a pair will be the $p$ points $\{ (x_j, y_j) \}$ in $\R^2_+$,
 where $\R_+$ denotes the positive reals.

The pair $(X^2,Y)$ will be in ${\frak C}_n^2$, the set of commuting pairs of Hermitian matrices.
Its spectrum will consist of the joint eigenvalues $\{ (x_j, \pm y_j) : 1 \leq j \leq p \} \subset \R_+ \times \R$.
 
 \section{Distribution of the skew spectrum}
 \label{secd}
 
 Let $n = 2p$.
 Let $(\R^2_+)^p_{\rm gen}$ be the set of $p$-tuples $\{ (x_j, y_j ) : 1 \leq j \leq p \}$ in $\R^2_+$ such
 that
all the $x_j$'s are distinct, and all the $y_j$'s are distinct.
Let $\ang$ denote the generic elements of $\an$, as described in Section \ref{secc}.
The map from $\un \times \rtpg$ to $\ang$ that sends $(U, \{ (x_j, y_j ) \})$ to
\[
(U \left[ \oplus_{j=1}^p  \begin{pmatrix}
  x_j & 0 \\
  0 & -x_j
  \end{pmatrix} \right] U^*,
  U \left[ \oplus_{j=1}^p  \begin{pmatrix}
  0 & y_j \\
  y_j & 0
  \end{pmatrix} \right] U^*)
  \]
  will not be injective, as any unitary $U$ that is the direct sum of $2$-by-$2$'s of the form  $ \begin{pmatrix}
  e^{ i \theta_j}  & 0 \\
  0 & e^{i \theta_j}
  \end{pmatrix} $ will leave the block-diagonal matrices invariant.
  To rectify this, we shall consider
  \begin{eqnarray}
  \notag
  G: \un / \T^p \times \rtpg &\to & \ang \\
  (U, \{ (x_j, y_j ) \}) & \mapsto & 
  (U \left[ \oplus_{j=1}^p  \begin{pmatrix}
  x_j & 0 \\
  0 & -x_j
  \end{pmatrix} \right] U^*,
  U \left[ \oplus_{j=1}^p  \begin{pmatrix}
  0 & y_j \\
  y_j & 0
  \end{pmatrix} \right] U^*).
  \label{eqd1}
  \end{eqnarray}
  Then $G$ is a bijection, and it follows from the proof of Theorem \ref{thmd1} that it is a diffeomorphism as
  the Jacobian does not vanish.
Let $\nu$ be volume measure on the homogeneous space $\un / \T^p$.
  
  To reduce the use of superscripts, we shall let $Z$ be an element of $\an$, and write its two components
  as
  \[
  Z = (Z^1, Z^2)  \= (X, Y) .
  \]
  If $x = \{x_1, \dots, x_p\} \in \C^p$, let 
  \beq
  A_x &\=&  \oplus_{j=1}^p  \begin{pmatrix}
  x_j & 0 \\
  0 & -x_j
  \end{pmatrix}  \\
  B_x  &=& \oplus_{j=1}^p  \begin{pmatrix}
 0&  x_j  \\
  x_j & 0
  \end{pmatrix} .
\eeq
  \bt
\label{thmd1}
Let $Z =  (X,Y)$ be chosen randomly in $\an$ with distribution $w(\| Z\|_F)$. Then the probability 
distribution of the skew spectrum of $(X,Y)$ on $\R^{2p}_+$  is given by
\begin{multline}
\label{eqd2}
\rho_n (x_1, y_1, \dots, x_p, y_p) \= \\
C_n\,  w(\| Z\|_F)\ \prod_{1 \leq k \leq p} \,  x_k y_k\sqrt{x_k^2 + y_k^2} 
\prod_{1 \leq i < j \leq p} 
\left[(x_i - x_j)^2 +(y_i-y_j)^2 \right]
\left[(x_i + x_j)^2 +(y_i-y_j)^2 \right] \\
\times
\left[(x_i - x_j)^2 +(y_i + y_j)^2 \right]
\left[(x_i + x_j)^2 +(y_i +y_j)^2 \right].
\end{multline}
\et
  \bp
Fix some point $Z = (X,Y) \in \ang$. Without loss of generality, we can choose a basis so that
$U$ in \eqref{eqd1} is the identity, and $(X,Y)$ has the form \eqref{eqa1} with skew spectrum in $\rtpg$.
 The derivative of $G$ is  a map
between the tangent spaces.
\[
dG : ( T_{[I \T^p]} \  \un / \T^p)  \times \R^{2p}  \ \to \ T_{(X,Y)} \an .
\]
If we view $dG$ as a real linear map, then the Jacobian will be ${\mathcal J} = \sqrt{\det (dG^* dG)}$.
So we will have
\be
\label{eqd3}
w(\| Z \|_F) dZ \=  w(\| Z \|_F) {\mathcal J}\ d\nu dx_1 dy_1 \dots dx_p dy_p .
\ee
Integrating with respect to $\nu$, we get that
$\rho_n$ equals the Jacobian times $w(\| Z \|_F)$; we must prove that this has the form \eqref{eqd2}.

The tangent space at the identity of $\un$ is the space of  skew-symmetric matrices. The tangent space of
$\un /\T^p$ at $[I \T^p]$ is the skew-symmetric matrices whose diagonals are of the form $(\pm \ii \theta_j )$.
(We shall write $\ii$ for $\sqrt{-1}$ to distinguish from $i$ used as an index).
We have
\begin{eqnarray}
\label{eqd4}
\notag
dG|_{(I \T^p, x,y) } (S, a, b) &\=& \frac{d}{dt} \left( e^{tS} A_{x+ta} e^{-tS}, e^{tS} B_{y+tb} e^{-tS} \right) \\
&=& \left( SA_x - A_x S  + A_a, SB_y - B_y S + B_b \right).
\end{eqnarray}

We want to pick a  basis for the tangent space that facilitates computation. 

For any $k,\ell \in \{ 1, \dots, n\}$ let $E_{k,l}$ denote the elementary matrix with $1$ in the $(k,\ell)$ place and $0$ elsewhere.
We shall let $i,j$ range between $1$ and $p$, and $\alpha$ and $\beta$ range over $\Z_2$ (where $1+1 = 0$).
Define a basis as follows. 

For each $1 \leq k \leq p$, we have 3 matrices:
\beq
R_k &\= & \frac{1}{\sqrt{2}} [ E_{2k-1,2k} - E_{2k,2k-1} ]\\
S_k &=& \frac{\ii}{\sqrt{2}} [ E_{2k-1, 2k} + E_{2k,2k-1}]\\
T_k &=& \frac{\ii}{\sqrt{2}} [ E_{2k-1, 2k-1} - E_{2k,2k}] .
\eeq
 For each pair $(i,j)$ in $\{ 1, \dots , p \}$ with $i < j$ and each $\alpha, \beta \in \Z_2$, we have two matrices
 \beq
 R_{ij,\alpha\beta} &\=& \frac{1}{\sqrt{2}} [ E_{2i-\alpha, 2j - \beta} - E_{2j-\beta, 2i-\alpha}]\\
  S_{ij,\alpha\beta} &\=& \frac{\ii}{\sqrt{2}} [ E_{2i-\alpha, 2j - \beta} - E_{2j-\beta, 2i-\alpha}].
  \eeq
Finally, for a basis of $\R^{2p}$, thought of as the tangent space to $\rtpg$, we let $\{e_k^1: 1 \leq k \leq p \}$ be the standard basis for the first slot, and $\{e_k^2: 1 \leq k \leq p \}$ be the standard basis in the second slot.

 Straightforward calculations show
 \begin{eqnarray*}
 &[R_k, A_x] &= 2\ii x_k S_k,\qquad [R_k, B_y] = -2\ii y_k T_k  \\
&[S_k, A_x]  &= -2\ii x_k R_k, \quad [S_k, B_y] = 0 \\
&[T_k, A_x]  &= 0, \qquad\qquad [T_k, B_y] = 2\ii y_k R_k
\end{eqnarray*}
and
\beq
& [R_{ij,\alpha \beta}, A_x]&=  \ii ((-1)^\beta x_j - (-1)^\alpha x_i) S_{ij,\alpha \beta}  \\
& [ R_{ij,\alpha \beta}, B_y]&=  \ii (y_i S_{ij,\alpha+1,\beta} - y_j S_{ij,\alpha,\beta+1} )\\
& [S_{ij,\alpha \beta}, A_x]&=  \ii ((-1)^\alpha x_i   - (-1)^\beta x_j) R_{ij,\alpha \beta} \\
& [S_{ij,\alpha \beta}, B_y]&=  \ii ( y_j R_{ij,\alpha, \beta+1} - y_i R_{ij,\alpha+1, \beta} ).
 \end{eqnarray*}
The matrix for $dG$ has a block form. It maps the 5 (real) dimensional space spanned by  
$\{ R_k, S_k, T_k, e^1_k, e^2_k \}$ into the 6 dimensional space that is spanned by $\{ R_k, S_k, T_k \}$ in both the first and second slots,
and it maps the 8 dimensional space spanned by $\{  R_{ij,\alpha\beta},  S_{ij,\alpha\beta} : \alpha, \beta \in \Z_2 \}$ 
into two copies of the same space.
Because of the block structure, the Jacobian of the whole map will be the product of the Jacobians for each block.

The first set of blocks look like this.
\be
\label{eqd6}
\bordermatrix{~&R_k&S_k&T_k&e_k^1&e_k^2 \cr
R_k&0 &-2x_k &0 &0 & 0
\cr
S_k&2 x_k &0 &0 & 0& 0
\cr
T_k&0 &0 & 0& -\sqrt{2}& 0
\cr
R_k& 0&0 & 2y_k&0 & 0
\cr
S_k& 0&0 &0 &0 & -\sqrt{2}
\cr
T_k& -2  y_k &0 &0 &0 & 0
\cr} \ \times \ii
\ee
The second set looks like
\be
\label{eqd7}
\bordermatrix{~&R_{ij00}&R_{ij10}&R_{ij01}&R_{ij11}&S_{ij00}&S_{ij10}&S_{ij01}&S_{ij11} \cr
R_{ij00}&0&0 & 0&0 				&x_i-x_j &0 &0 &0 
\cr
R_{ij10}&0&0 & 0&0				 &0 &-x_i-x_j &0 &0 
\cr
R_{ij01}&0&0 & 0&0 				&0 &0 &x_i + x_j &0 
\cr
R_{ij11}&0&0 & 0&0				 &0 &0 &0 & -x_i + x_j
\cr
S_{ij00}&x_j - x_i&0&0 & 0            &0 &0 &0 &0 
\cr
S_{ij10}&0&x_j +x_i & 0&0           &0 &0 &0 &0 
\cr
S_{ij01}&0&0 &-x_j - x_i  & 0       &0 &0 &0 &0 
\cr
S_{ij11}&0&0&0&-x_j +x_i           & 0&0 &0 &0 
\cr
R_{ij00}&0&0 & 0&0 				&0 &-y_i &y_j &0 
\cr
R_{ij10}&0&0 & 0&0 				&-y_i &0 &0 &y_j 
\cr
R_{ij01}&0&0 & 0&0 				&y_j &0 &0 & -y_i
\cr
R_{ij11}&0&0 & 0&0 				&0 &y_j &-y_i &0 
\cr
S_{ij00}&0&y_i & -y_j &0 	&0 &0 &0 &0 
\cr
S_{ij10}&y_i &0&0 & -y_j 		&0 &0 &0 &0 
\cr
S_{ij01}&-y_j &0&0 &y_i		 &0 &0 &0 &0 
\cr
S_{ij11}&0&-y_j & y_i&0 		&0 &0 &0 &0 
\cr
} \ \times \ii
\ee
When \eqref{eqd6} is premultiplied by its adjoint, one gets a diagonal matrix whose determinant is 
\be
\label{eqd8}
\notag
4^4 (x_k^2 + y_k^2)x_k^2 y_k^2 .
\ee
When \eqref{eqd7} is premultiplied by its adjoint, the resulting matrix has determinant
\be
\label{eqd9}
\notag
\left[\left( (x_j - x_i)^2 +y_j^2 + y_i^2 \right)^2 - 4 y_j^2y_i^2 \right]^2
\left[ \left( (x_j + x_i)^2 +y_j^2 +y_i^2 \right)^2 - 4 y_j^2 y_i^2 \right]^2 ,
\ee
which equals the square of 
\[
\left[(x_i - x_j)^2 +(y_i-y_j)^2 \right]
\left[(x_i + x_j)^2 +(y_i-y_j)^2 \right] 
\left[(x_i - x_j)^2 +(y_i + y_j)^2 \right]
\left[(x_i + x_j)^2 +(y_i +y_j)^2 \right].
\]

Multiplying all these together, we get that the Jacobian times $w(\| Z \|_F)$, when integrated
with respect to $\nu$, is \eqref{eqd2}.
  \ep
  
%
  
  \begin{lemma}
  \label{lemd1}
  Define
  \begin{multline}
  f(x_i,x_j,y_i,y_j) \= 
\left[(x_i - x_j)^2 +(y_i-y_j)^2 \right]
\left[(x_i + x_j)^2 +(y_i-y_j)^2 \right] \\
\times
\left[(x_i - x_j)^2 +(y_i + y_j)^2 \right]
\left[(x_i + x_j)^2 +(y_i +y_j)^2 \right].
 \end{multline}
  Let $d^2 = (x_i-x_j)^2 + (y_i - y_j)^2$.
 Let $\vare, M$ be positive constants, and assume that $x_1,x_2,y_1$ and $y_2$ are all between
 $\vare$ and $M$. Then we have
 \be
 \label{eqd10}
 128 \vare^6 \ d^2 \ \leq\  f (x_i,x_j,y_i,y_j)\
 \leq \ 200 M^6\  d^2.
 \ee

 \end{lemma}
 \bp
The first factor of $f$ is $d^2$. The other three factors are bounded below by
$ (2 \vare)^2 (2 \vare)^2 (8 \vare^2)$ and bounded above by
$(5M^2)(5M^2)(8M^2)$.
\ep
  It follows from Lemma \ref{lemd1}, that in compact subsets of $(0,\i)^2$ 
  the repulsion between elements of the skew spectrum, as given by \eqref{eqd2},
  is of the order of the square of their Euclidean distance.

  \section{Fekete points and the limiting distribution}
  \label{sece}
  
\begin{figure}[!htpb]
\includegraphics
[width = 0.40\paperwidth]
{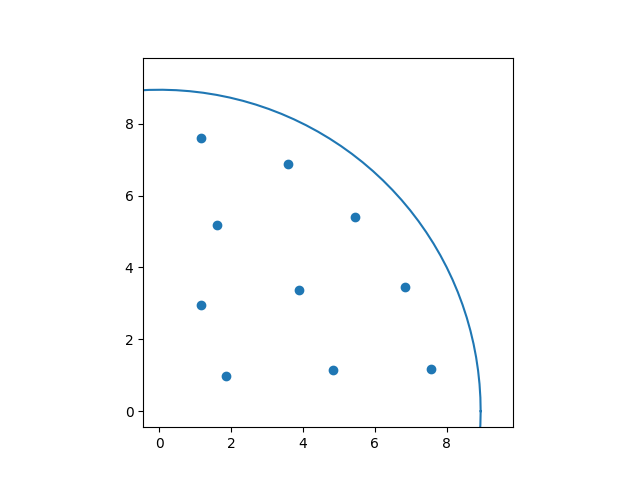}
\includegraphics
[width = 0.40\paperwidth]
{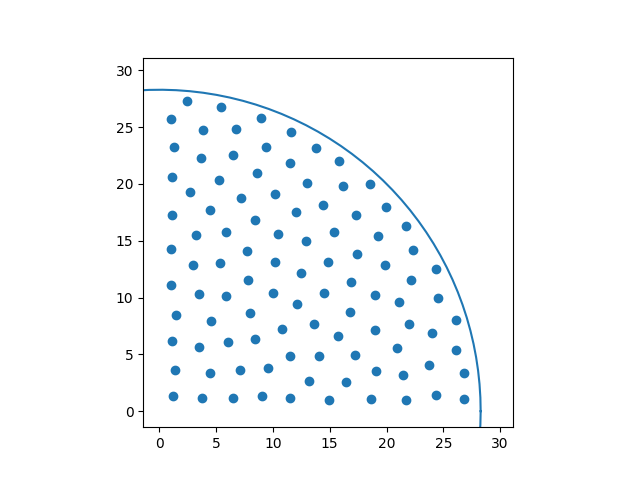}
\includegraphics
[width = 0.40\paperwidth]
{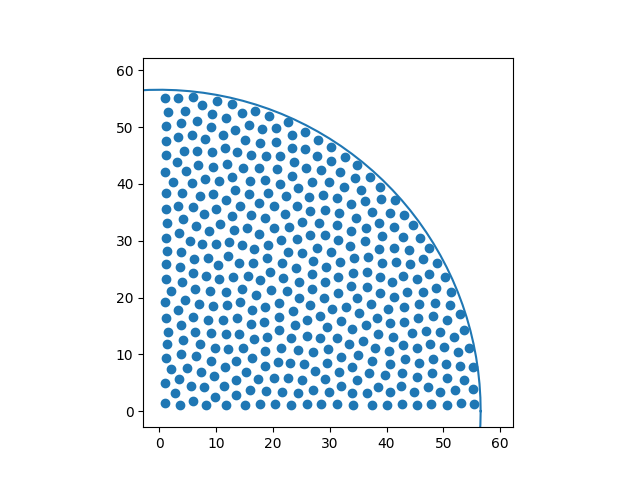}
\caption{Numerical simulation of points that maximize $\rho_n$ for  a pair of anti-commuting self-adjoint matrices,
with sizes $n=20, 200, 800$. Circle has radius  $2\sqrt{n}$.}
\label{fige1}
\end{figure}

 In this section, we shall just  consider the Gaussian case $w(Z) = e^{-\frac 12 \| Z\|_F^2}$.
 For any fixed $n=2p$, the skew spectrum is most likely to occur where the density $\rho_n$ from
 \eqref{eqd2} is highest. Using gradient descent, we numerically calculated what distribution of points maximized
 $ (\rho_n)$. See Figure \ref{fige1}.
They seem to be approximately equally spaced within the quarter-disk of radius $2\sqrt{n}$.

\begin{figure}[!htpb]
\includegraphics
[width = 0.40\paperwidth]
{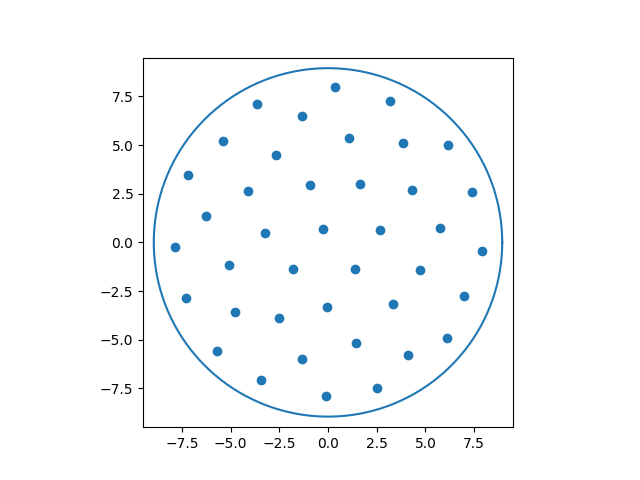}
\includegraphics
[width = 0.40\paperwidth]
{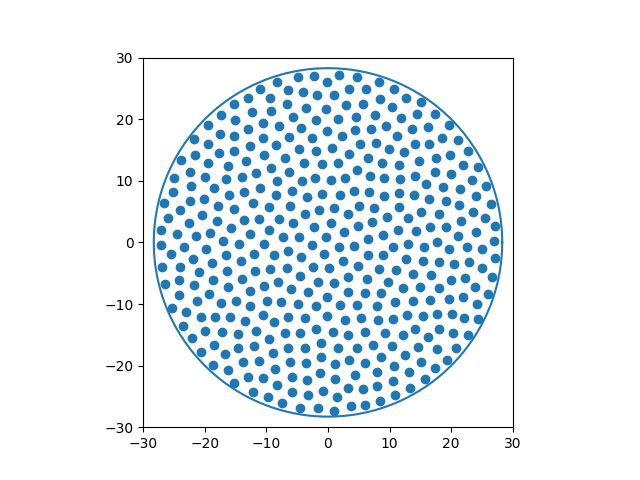}
\caption{Numerical simulation of points that maximize $\kappa_n$ for  a pair of commuting self-adjoint matrices,
with sizes $n=40, 400$. Circle has radius  $\sqrt{2n}$.}
\label{fige2}
\end{figure}

By way of comparison, we plot the joint eigenvalues that maximize the density $\kappa_n$ from \eqref{eqb2} for a pair
of commuting self-adjoint matrices, with the same Gaussian weight. In this case, they are approximately equally spaced within the disk of radius $\sqrt{2n}$. See Figure \ref{fige2}.

Write $z = (x,y)$ for a point in $\R_+^2$. 
\begin{definition}
Let $n = 2p$. We shall say a set $S = \{ z_1, \dots, z_p \} \subseteq \R_+^2$ is a maximal likelihood set 
if $\rho_n$ attains its maximum on $S$.
\end{definition}

It is not immediately obvious that maximal likelihood sets exist, since the domain is not bounded, but we shall prove that they do. It is more convenient to work with $\tau := - \log \rho_n$, a function
from $\R_+^{2p}$ to $(-\infty, \infty]$ which we want to minimize.  Let $f$ be as in Lemma \ref{lemd1}.
Then
\[
\tau(z_1, \dots , z_p) \=
\frac{1}{2} \sum_{k=1}^p |z_k|^2 - \log | x_k y_k z_k| 
- \frac{1}{2} \sum_{ k \neq \ell} \log f(z_k, z_\ell) .
\]

\begin{lemma}
\label{leme1}
For each $n \geq 1$, there exists a set $S$ so that $\tau(S) \leq n^2$.
\end{lemma}
\bp
Case: $p = q^2$ for $q \in {\mathbb N}$.
Let $S_q = \{ 1, 2, \dots, q \} \times \{1, 2, \dots, q \}$. By Lemma \ref{lemd1}, for any $k \neq \ell$
we have $f(z_k, z_\ell) \geq 1$. 
So 
\beq
\tau (S_q) & \ \leq \ & \frac{1}{2} \sum_{i,j=1}^q (i^2 + j^2) \\
&=& \frac{q^2(q+1)(2q+1)}{6} \\
&\leq & q^4 = \frac 14 n^2.
\eeq

General case: choose $q$ so that $(q-1)^2 < p \leq q^2$. Let $S$ be any $p$ elements of $S_q$.
Then
\beq
\tau(S) &\ \leq \ & q^4 \\
&\leq& 4 ((q-1)^2 +1)^2 \\
&\leq & 4 p^2 = n^2 .
\eeq
\ep

\begin{lemma}
\label{leme2}
Let $n = 2p$ be a positive even integer.
Choose $K \geq 3 p$ so that 
\be
\label{eqe1}
\frac{1}{2} K^2  - (3 p + 4p^2) \log ( K ) - \frac{p^2}{2} \log (400) - 4p^2 > 0.
\ee
Let $S  = \{ z_1, \dots, z_p \} \in (\R_+^2)^p$ be a set so that $\tau(S) \leq 4 p^2$.
Then the maximum length of an element of $S$ is at most $K $.
\end{lemma}
\bp
Let $M = \max \{ |z_k| : 1 \leq k \leq p \}$.
By Lemma \ref{lemd1}, for $k \neq \ell$ we have
\[
f(z_k, z_\ell) \ \leq \ 400 M^8.
\]
So
\beq
4p^2 &\ \geq \ & 
\tau(S)\\
 &=& \frac{1}{2} \sum_{k=1}^p \left(  |z_k|^2 - \log | x_k y_k z_k| \right) 
- \frac{1}{2} \sum_{ k \neq \ell} \log f(z_k, z_\ell) \\
 &\ \geq \ & \frac 12 M^2 - p \log M^3 - \frac{p(p-1)}{2} \log (400 M^8) \\
 &\geq &
 \frac 12 M^2 -(3p+ 4p^2) \log M - \frac{p^2}{2} \log 400.
\eeq
The last inequality
\be
\label{eqe2}
4p^2 \ \geq \ \frac 12 M^2 -(3p+ 4p^2) \log M - \frac{p^2}{2} \log 400
\ee
fails at $M = K$ (by choice of $K$). Moreover, the right-hand side of \eqref{eqd2}
is increasing for $M \geq 3p$, so we must have $M < K$.
\ep

\bt
\label{thme1}
Maximal likelihood sets exist.
\et
\bp
Let $K$ be as in Lemma \ref{leme2}.
Consider $\tau : [0,K]^{2p} \to (-\i,\i]$.
This is a continuous function on a compact set, so attains its infimum. 
By Lemmas \ref{leme1} and \ref{leme2} this is a global infimum for $\tau$ on $\R_+^{2p}$.
\ep

\begin{definition}
A set $S \subseteq (\R_+^{2})^p$ is a Fekete set of size $p$ if the set $\sqrt{p} S$ is a maximal likelihood 
set of $\rho_n$.

A measure is a Fekete measure of size $p$ if it consists of $p$ atoms of weight $\frac 1p$ at each point of
a Fekete set of size $p$.
\end{definition}

Based on the numerical simulations shown in Figure \ref{fige1}
and analogy with the commuting case, we are led to ask the following questions.

\begin{question}
\label{qe4}
Let $\mu_p$ be a sequence of Fekete measures of size $p$. 
\begin{enumerate}
\item
Is there a compact set that contains the support of every $\mu_p$?
\item
Does the sequence
$\mu_p$ converge weakly (when integrated against bounded continuous functions) to a unique
compactly supported measure $\mu$?
 \item
 Define a random probability measure $\nu_p$ by choosing random anti-commuting self-adjoint pairs of size $2p$-by-$2p$ with the Gaussian measure  $C_p e^{-\frac 12 \| Z \|^2} dZ$, and letting $\nu_p$ 
 have mass $\frac 1p$ at each point of the skew spectrum  of  $Z$. Let ${\mathbb E}_p$ denote expectation with respect to this process.  Is there a measure $\mu$ so that 
  \[
  \lim_{p \to \i} {\mathbb E}_p \int_{\R^2_+}  \phi(z) d\nu_p(z)  \= \int_{\R^2_+} \phi(z) d \mu(z)  \qquad \forall \phi \in C_b(\R^2_+)?
  \]
  \item
If the answer to Questions 2 and 3 is yes, is  $\mu$ equal to normalized  area measure on the quarter disk in the first quadrant of radius $\sqrt{8}$?
\end{enumerate}
\end{question}

\end{document}

%% file: def_latex5.tex
\def\mcc{M\raise.5ex\hbox{c}C}
\def\mccarthy{M\raise.5ex\hbox{c}Carthy}

\def\l{\lambda}

\def\vare{\varepsilon}

\let\ii=\i
\let\i=\infty

\def\={\ = \ }

\def\C{\mathbb C}
\def\R{\mathbb R}
\def\T{\mathbb T}

\def\be{\begin{equation}}
\def\ee{\end{equation}}
\def\beq{\begin{eqnarray*}}
\def\eeq{\end{eqnarray*}}

\def\bp{\begin{proof}}
\def\ep{\end{proof}}

\def\bl{\begin{lemma}}
\def\el{\end{lemma}}
\def\bt{\begin{theorem}}
\def\et{\end{theorem}}
\def\bprop{\begin{proposition}}
\def\eprop{\end{proposition}}
\def\bd{\begin{definition}}
\def\ed{\end{definition}}
\def\br{\begin{remark}}
\def\er{\end{remark}}
\def\bexer{\begin{exercise}}
\def\eexer{\end{exercise}}

\newtheorem{theorem}{Theorem}[section]
\newtheorem{proposition}[theorem]{Proposition}
\newtheorem{lemma}[theorem]{Lemma}

\newtheorem{remark}[theorem]{Remark}

\newtheorem{definition}[theorem]{Definition}

\newtheorem{question}[theorem]{Question}